\documentclass[11pt]{amsart}

\usepackage{amsmath}
\usepackage{amssymb}
\usepackage{graphicx}

\newtheorem{theorem}{Theorem}[section]
\newtheorem{corollary}[theorem]{Corollary}

\newtheorem{proposition}[theorem]{Proposition}
\theoremstyle{definition}
\newtheorem{definition}[theorem]{Definition}
\newtheorem{example}[theorem]{Example}

\numberwithin{equation}{section}

\newcommand \s{^{*}}
\newcommand \+{^{\dag}}
\newcommand \p{^{\perp}}

\title[Complementable Operators and their Schur Complements]{Complementable Operators and their Schur Complements}

\author[Sachin Manjunath Naik]{Sachin Manjunath Naik}

\address[Sachin Manjunath Naik]{Department of Mathematical and Computational Sciences, National Institute of Technology Karnataka (NITK), Surathkal, Mangaluru 575 025, India}
\email{{\tt sachinmaths46@gmail.com}}

\author[P. Sam Johnson]{ P. Sam Johnson}
\address[P. Sam Johnson]{Department of Mathematical and Computational Sciences, National Institute of Technology Karnataka (NITK), Surathkal, Mangaluru 575 025, India}
\email{\tt sam@nitk.edu.in}



\begin{document}

\begin{abstract}
In this paper, we characterize complementable operators and provide more precise expressions for the Schur complement of these operators using a single Douglas solution. We demonstrate the existence of subspaces where the given operator is invariably complementable. Additionally, we investigate the range-Hermitian property of these operators.   
\end{abstract}

\maketitle

\begin{sloppypar}
\section{Introduction}

Let $T$ be an $(n+m)\times (n+m)$ matrix over an arbitrary field which is partitioned in the form
$$T=\begin{pmatrix} A & B \\
C & D 
\end{pmatrix}$$
with a nonsingular matrix $D$ of order $m$.  Then, the classical Schur complement of $D$ in $T$ is the $n\times n$ matrix defined by 
\begin{equation}\label{equation-1}
T_{/D} = A-BD^{-1}C.
\end{equation}
The formula (\ref{equation-1}) was first used by Schur \cite{Schur}.  However, this terminology is due to Haynsworth \cite{Emilie}.  Matrices of the form (\ref{equation-1}) have undoubtedly been encountered from the time matrices were first used, particularly in Gaussian elimination \cite{Gantmacher, Zhang} and minor determinants \cite{sylvester} and others. We remark that the Schur complement in this classical sense is a matrix of lower order, while the Schur complement in our sense is always a matrix of the order of $T$ and it is henceforth denoted by 
$$T_{/A}=\begin{pmatrix} A-BD^{-1}C & 0 \\
0 & 0 
\end{pmatrix}.$$
The Schur complement has been extended to the case where the block matrix $D$ is rectangular and singular, by replacing the regular inverse by the Moore-Penrose inverse or any other generalized inverse\cite{Albert1}.   This notion has proved to be a fundamental idea in many applications like numerical analysis, statistics, operator inequalities, electrical network theory, discrete-time regulator problem and to name a few.  

Let $\mathcal H$ and $\mathcal K $ be real or complex Hilbert spaces.  We denote  $\mathcal B(\mathcal H, \mathcal K)$  the space of all linear bounded operators from $\mathcal H$ to $\mathcal K$ and we abbreviate $\mathcal B(\mathcal H)=\mathcal B(\mathcal H, \mathcal H)$. For $T \in \mathcal B(\mathcal H, \mathcal K),$ we denote by $T^*$, $\mathcal N(T)$ and $\mathcal R(T),$ respectively, the adjoint, the null-space and the range-space of $T$. An operator $T\in \mathcal B(\mathcal H)$ is said to be a projection if $T^2=T.$ A projection, $T$ is said to be orthogonal if $T=T^*.$ An operator $T\in \mathcal B(\mathcal H)$ is called positive (denoted $T\geq 0$) if $T=T^*$ and  $\langle Tx, x \rangle \geq 0$ for all $x\in \mathcal H$.  The cone of positive operators in $\mathcal B(\mathcal H)$ is denoted by $\mathcal B(\mathcal H)^+$.  The operator $|T|:=(T^*T)^{1/2}$ is called the modulus of $T$. 

Anderson \cite{Anderson} first introduced the definition of the Schur complement for positive operators on Hilbert spaces. Given $T\in \mathcal B(\mathcal H)^+$, the shorted operator (or Schur complement) of $T$ with respect to a closed subspace $M$ of $\mathcal H$ is defined as  
$$T_{/M} =\max _{\leq }\{X \in \mathcal B (\mathcal{H})^+: X \leq T~ \text{and} ~ \mathcal{R}(X) \subseteq M^\perp \}.$$ The existence of maximum has been studied by Krein \cite{Krein} in the study of symmetric operator extensions. Ando \cite{Ando} introduced the definition of the Schur complement for $M$-complementable bounded operators.  Since then, the Schur complement of bounded linear operators has been attracted extensively by many researchers. See, for example \cite{Carlson1,Corach2,Corach1,Giribet, Pekarev,Ptak}.
Ando's definition opened the door to various extensions and generalizations, as observed by Butler and Morley \cite{Butler}, who compared six different notions of Schur complements for finite dimensions. 

Mitra \cite{Mitra} generalized the concept of complementable operators to the context of operators between distinct Hilbert spaces, introducing bilateral shorted operators and exploring their properties. Subsequently, Antezana et al. \cite{Antezana} formulated the concept of complementability and weak complementability for bounded operators $T:\mathcal{H} \to \mathcal{K}$ relative to closed subspaces $M \subseteq \mathcal{H}$, $N \subseteq \mathcal{K}$, which gave rise to the bilateral shorted operator and its applications in extending parallel sums of bounded linear operators. A comprehensive survey of the Schur complement can be found in \cite{Arias,Carlson1,Cottle,Ouellette}.   

We now describe the contents of the paper.  In the next section, we provide a brief background for the rest of the material in the paper. In Section 3, we comprehend complementable bounded operators with more characterizations. Furthermore, for complementable operators relative to closed subspaces $M$ and $N$, we present a new approach to derive the Schur complement using a single Douglas solution. In Section 4, we show that this methodology provides insights into the closed range of complementable operators and their Schur complements relative to the associated subspaces. Also, we provide an extension of a characterization of $EP$-matrix and its Schur complement given in \cite{Meenakshi} to a hypo-$EP$ operator $T:\mathcal{H} \to \mathcal{H}$ relative to closed subspaces $M \subseteq \mathcal{H}$ and  $N \subseteq \mathcal{H}$.

\section{Preliminaries}
	
Let $M$ and $N$ be  closed subspaces of the Hilbert spaces $\mathcal{H}$ and $\mathcal{K}$, respectively.  A bounded linear operator $T\in \mathcal B(\mathcal H, \mathcal K)$ can be expressed as 
$$T=\begin{pmatrix}
A & B\\
C & D
\end{pmatrix}$$
with regard to the  decompositions  $\mathcal{H}=M \oplus M\p$ and $\mathcal{K}=N \oplus N\p$.

\begin{definition}\cite{Antezana}
	Let $P_r \in \mathcal{B}(\mathcal{H})$ and $P_\ell \in \mathcal{B}(\mathcal{K})$ be projections. An operator $T \in \mathcal {B}(\mathcal{H}, \mathcal{K})$ is called $(P_r, P_\ell)$-complementable if there exist operators $M_r \in \mathcal{B}(\mathcal{H})$ and $M_\ell \in \mathcal{B}(\mathcal{K})$ such that 
	\begin{enumerate}
		\item $(I-P_r)M_r=M_r$ and $(I-P_\ell)TM_r=(I-P_\ell) T$,
		\item $M_\ell (I-P_\ell)=M_\ell$ and $M_\ell T(I-P_r)=T(I-P_r) $.
	\end{enumerate}
\end{definition}

\begin{proposition}\cite{Antezana} \label{propdef}
	Let $P_r \in \mathcal{B}(\mathcal{H})$ and $P_\ell \in \mathcal{B}(\mathcal{K})$ be  projections with their ranges $M$ and $N,$ respectively.
	Then $T$ is $(P_r, P_\ell)$-complementable if and only if $\mathcal{R}(C)\subseteq \mathcal{R}(D)$ and $\mathcal{R}(B\s) \subseteq \mathcal{R}(D\s)$.
\end{proposition}

It is observed from Proposition \ref{propdef} that $(P_r, P_\ell)$-complementable depends only on the ranges of $P_r $ and $ P_\ell$.

\begin{definition}\cite{Antezana}
	An operator $T \in \mathcal {B}(\mathcal{H}, \mathcal{K})$ is called $(M, N)$-complementable if it is $(P_r, P_\ell)$-complementable for some projections $P_r $ and $ P_\ell$ with $\mathcal{R}(P_r)=M$ and $\mathcal{R}(P_\ell)=N$.
\end{definition}

\begin{definition}\cite{Antezana} An operator $T \in \mathcal {B}(\mathcal{H}, \mathcal{K})$ is called $(M,N)$-weakly complementable if
	 $$\mathcal{R}(C)\subseteq \mathcal{R}(|D\s|^{\frac{1}{2}}) \ \text{ and } \ \mathcal{R}(B\s) \subseteq \mathcal{R}(|D|^{\frac{1}{2}}).$$
\end{definition}
We know that $\mathcal{R}(S)\subseteq \mathcal{R}(|S^*|^{1/2})$, for any bounded operator $S$.  By Proposition \ref{propdef}, it is easy to see that if $T$ is  $(M,N)$-complementable, then it is $(M,N)$-weakly complementable. Both notions of complementability coincide if $\mathcal {R}(D)$ is closed from the following facts : 
\begin{enumerate}
	\item[(i)] $\mathcal R(D)$ is closed if and only if $\mathcal R(D^*)$ is closed ;
	\item[(ii)] $\mathcal R(D)=\mathcal R(D^{1/2})$ if and only if $\mathcal R(D)$  is closed ;
	\item[(iii)] $\mathcal R(D)=\mathcal R(|D^*|)= \mathcal R(|D^*|^{1/2})$.
\end{enumerate}
From the above observations, it is noted that if 
$T=\begin{pmatrix}
A & B\\
C & D
\end{pmatrix}$ is $(M,N)$-complementable with $\mathcal R(D)$ closed, then 
$$T_{/(M,N)}=\begin{pmatrix}
A-BD\+ C & 0\\
0 & 0
\end{pmatrix}.$$

\begin{theorem}\label{dgls}
	Let $A, B \in \mathcal{B}(\mathcal H)$.  Then the following are equivalent:
	\begin{enumerate}
		\item $\mathcal{R}(A) \subseteq \mathcal{R}(B)$ ;
		\item $AA\s \leq \lambda BB\s$, for some $\lambda >0$ ;
		\item There exists a bounded operator $C \in \mathcal{B}(\mathcal H)$ such that $A=BC$.
	\end{enumerate}
	Moreover, if these equivalent conditions hold, then there is a unique operator $C \in \mathcal{B}(\mathcal H)$ such that
	\begin{enumerate}
		\item[(i)] $\|C\| = \inf \{ \lambda >0: AA\s \leq \lambda BB\s\} $;
		\item[(ii)] $\mathcal{N}(A) = \mathcal{N}(B)$;
		\item[(iii)] $\mathcal{R}(C) \subseteq \mathcal{N}(B)\p$.
	\end{enumerate}
\end{theorem}

\begin{definition}\cite{Antezana}
Let $T$  be $(M,N)$-weakly complementable. Let $F$ and $E$ be the reduced solutions of the equations $C=|D\s|^{\frac{1}{2}}UX$ and $B\s= |D|^{\frac{1}{2}}X$, respectively, where $D=U|D|$ is the polar decomposition of $D$. Then the bilateral shorted operator (the Schur complement) of $T$ with respect to closed subspaces 
	$M$ and $N$ is   $ \begin{pmatrix}
		A-E\s F & 0 \\
		0 & 0
	\end{pmatrix}$ and is denoted by $T_{/(M,N)}$.\\ If $\mathcal{H}=\mathcal K$ and $M=N$, then we denote $T_{/(M,N)}$ by $T_{/M}$. 
\end{definition}

The above definition relies on the existence of two reduced solutions to define the Schur complement of weakly complementable operators.  But for an 
$(M,N)$-complementable operator $T=\begin{pmatrix}
A & B\\
C & D
\end{pmatrix}$, the Schur complement can be elegantly represented by employing either one of the Douglas solutions pertaining to the conditions $\mathcal{R}(C)\subseteq \mathcal{R}(D)$ or $\mathcal{R}(B\s) \subseteq \mathcal{R}(D\s)$.

Consider the reduced Douglas solution $Z$
for the inclusion $\mathcal{R}(C)\subseteq \mathcal{R}(D)$.  As a result, we can express 
$C=DZ=U|D|Z=U{|D|^{\frac{1}{2}}}{|D|^{\frac{1}{2}}}Z={|D\s|^{\frac{1}{2}}}U{|D|^{\frac{1}{2}}}Z$. Therefore,  ${|D|^{\frac{1}{2}}}Z$ serves as a solution for the equation 
$C=|D\s|^{\frac{1}{2}}UX$. It is also important to note that $\mathcal{R}({|D|^{\frac{1}{2}}}Z) \subseteq \overline{\mathcal{R}({|D|^{\frac{1}{2}}})}$.

Therefore, ${|D|^{\frac{1}{2}}}Z$ stands as the reduced solution for the equation $C=|D\s|^{\frac{1}{2}}UX$. Furthermore, let $B=E\s |D|^{\frac{1}{2}}$. This allows us to express $BZ=E\s |D|^{\frac{1}{2}} Z=E\s F$ and hence $$T_{/(M,N)}=\begin{pmatrix}
	A-BZ & 0 \\
	0 & 0
\end{pmatrix}.$$ Moreover, if $Y$ is a solution for the equation $B=XD$, then $BZ=YC$.  Hence, the Schur complement can also be represented as $$T_{/(M,N)}=\begin{pmatrix}
	A-YC & 0 \\
	0 & 0
\end{pmatrix}.$$

In this context, the condition $\mathcal{R}(C) \subseteq \mathcal{R}(D)$ is essential for establishing the existence of the operator $Z$  such that $C=DZ$. Furthermore, through the following example, we can observe that the condition  $\mathcal{R}(B\s) \subseteq \mathcal{R}(D\s)$ ensures the Schur complement is well-defined.
Indeed, if we consider two solutions,  $Z_1$ and $Z_2$, such that $C=DZ_1$ and $C=DZ_2$, we find that $D(Z_1-Z_2)=0 $ if and only if $ \mathcal{R}(Z_1-Z_2) \subseteq \mathcal{N}(D) \subseteq \mathcal{N}(B) $, which implies $ BZ_1=BZ_2$. Therefore, it becomes evident that the definition of the Schur complement is independent of the specific choice of $Z$.

		\begin{theorem}\cite{Antezana}\label{C}
				Let $T$ be $(M,N)$-complementable. Then the following statements hold good :
				\begin{enumerate}
					\item $T\s$ be $(N,M)$-complementable ;
					\item $\mathcal{R}(T_{/(M,N)})=\mathcal{R}(T) \cap N$;
					\item $\mathcal{N}(T_{/(M,N)})= M\p + \mathcal{N}(T)$;
\item 	If $\mathcal{R}(D)$ is closed, then $T_{/(M,N)}=\begin{pmatrix}
					A-BD\+ C & 0\\
					0 & 0
		\end{pmatrix}.$
						\end{enumerate}
	\end{theorem}

	\begin{theorem}\cite{MR3314341} \label{t215}
		Let $T$ be $(M,N)$-complementable with $\mathcal{R}(C\s) \subseteq \mathcal{R}((T_{/(M,N)})\s)$, $\mathcal{R}(B) \subseteq \mathcal{R}(T_{/(M,N)})$ and $\mathcal{R}(T)$ closed. Then $$T\+=\begin{pmatrix}
			{T_{/M}}\+ & -{T_{/M}}\+ BD\+\\
			-D\+C{T_{/M}}\+ & D\+ +D\+ C{T_{/M}}\+BD\+
		\end{pmatrix}.$$
	\end{theorem}

\noindent 
For a complementable operator, the Schur complement is established through the Douglas solution of the inclusion $\mathcal{R}(C) \subseteq \mathcal{R}(D)$. Similarly, in the context of an operator satisfying $\mathcal{R}(C) \subseteq \mathcal{R}(D)$ and $\mathcal{R}(B^*) \nsubseteq \mathcal{R}(D^*)$, a Schur complement expression can be formulated, albeit not uniquely. 
We conclude this section by an illustrative example demonstrating these assertions.
	\begin{example}\label{ex5}
		Let $T=\begin{pmatrix}
		A & B\\
		C & D
		\end{pmatrix}:  \ell_2 \oplus \ell_2 \to \ell_2 \oplus \ell_2$ where $A$, $B$, $C$ and $D$ are defined from $\ell_2$ to $\ell_2$ as follows : 
		\begin{eqnarray*}
			A(x_1,x_2,\ldots) &= &\Big(0,x_1,\frac{x_2}{2},\ldots \Big);\\
			B(x_1,x_2,\ldots)&=&\Big(x_1,x_2,x_3,\ldots\Big);\\
			C(x_1,x_2,\ldots)&=&\Big(x_1,\frac{x_2}{2},\frac{x_3}{3}\ldots\Big);\\
			D(x_1,x_2,\ldots)&=&\big(x_2,x_3,x_4,\ldots\big).
		\end{eqnarray*}

		It is evident that $\mathcal{R}(D\s)=\ell_2\backslash\text{span}\{e_1\}$, $\mathcal{R}(C) \subseteq \mathcal{R}(D)=\ell_2$ and $\mathcal{R}(B\s)\nsubseteq \mathcal{R}(D\s)$. Define the operators $Z_1(x_1,x_2,\ldots)=\big(0,x_1,\frac{x_2}{2},\ldots\big)$ and $Z_2(x_1,x_2,\ldots)=\big(x_1,x_1,\frac{x_2}{2},\ldots\big)$. It is evident that $C=DZ_1=DZ_2$.\\ 
		Now, consider $$({T_{/M}}_1)(x_1,x_2,x_3,\ldots)=(A-BZ_1)(x_1,x_2,x_3,\ldots)=(x_1,0,0,\ldots)$$ and $$({T_{/M}}_2)(x_1,x_2,x_3,\ldots)=(A-BZ_2)(x_1,x_2,x_3,\ldots)=(0,0,0,\ldots).$$
		
		 In accordance with the definition of the Schur complement, both $({T_{/M}})_1$ and $({T_{/M}})_2$ qualify as Schur complements; however, it is important to note that $({T_{/M}})_1$ is not equal to $({T_{/M}})_2$.
	\end{example}

	\section{\textbf{Complementable Operators and Schur Complements}}
 Butler and Moorley's \cite{Butler} pioneering work illuminates a geometric pathway to comprehend complementable operators within finite dimensions. Drawing inspiration from their groundbreaking insights, our research embarks on an exploration into the vast domain of infinite-dimensional operators. By pushing the boundaries of traditional understanding, we aim to unveil new perspectives and insights that can reshape the landscape of operator theory and its applications.
\begin{theorem}\label{eqgm}
	
	Let $T \in \mathcal{B}(\mathcal{H},\mathcal{K})$ with $\mathcal{R}(D)$ closed. Then $T$ is $(M,N)$-complementable if and only if for each $x \in M$, there exists unique $z\in M$ such that $\{T(x,0)+T(M\p)\} \cap N = \{(z,0)\}$. Moreover, $T_{/(M,N)}(x,0)=(z,0)$.
\end{theorem}

\begin{proof}
	Let $Z$ be such that $C=DZ$. Then $T(x,0)+T(0,-Zx)=((A-BZ)x,0)$. Therefore, $\{T(x,0)+T(M^\perp)\} \cap N$ is non-empty.
	
	Now, assume $(z_1,0), (z_2,0) \in \{T(x,0)+T(M^\perp)\} \cap N$. There exist $y_1, y_2 \in M^\perp$ such that $T(x,y_1)=(z_1,0)$ and $T(x,y_2)=(z_2,0)$. Hence $Ax+By_1=z_1$, $Cx+Dy_1=0$, $Ax+By_2=z_2$, and $Cx+Dy_2=0$. Thus, $y_1-y_2 \in \mathcal{N}(D) \subseteq \mathcal{N}(B)$, leading to $By_1=By_2$. This implies $(z_1,0)=(z_2,0)$.
	
	Conversely, assume that $\mathcal{R}(D)$ is closed. For each $x\in M$, let $\{T(x,0)+T(M^\perp)\} \cap N = \{(z, 0)\}$ for some $z \in M$. Then there exists $y \in M^\perp$ such that $T(x,y)=(z,0)$, implying $Cx+Dy=0$. Thus $\mathcal{R}(C) \subseteq \mathcal{R}(D)$.
	
	Let $y\in \mathcal N(D)$. Then $T(0,0)+T(0,y)=(By,0) \in N$. Since $T(0,0)+T(0,0)=(0,0) \in N$ and $\{T(x,0)+T(M^\perp)\} \cap N$ is a singleton, we have $By=0$. Therefore $\mathcal{N}(D) \subseteq \mathcal{N}(B)$, leading to $\overline{\mathcal{R}(B^\ast)} \subseteq \overline{\mathcal{R}(D^\ast)}$. Since $\mathcal{R}(D)$ is closed, we obtain ${\mathcal{R}(B^\ast)} \subseteq {\mathcal{R}(D^\ast)}$.\\
	As observed earlier, we have ${T(x,0)+T(M^\perp)} \cap N = {z}=((A-BZ)x,0)$. Thus $T_{/(M,N)}(x,0)=(z,0)$.
\end{proof}
In the following example, we find Schur complement of an operator using our characterization.
\begin{example}
	Let $T:\ell_2 \rightarrow \ell_2$ be defined by
	\[ T(x_1,x_2,\ldots) = \Big(x_3+x_2-x_1,x_1+x_2,x_5+x_4-x_3,\frac{x_3}{3}+x_4,\ldots\Big) .\]
Consider the subspaces  $$M=\{x=(x_1,0,x_3,0,\ldots)~:~x \in \ell_2\}$$ and $$M^\perp =\{x=(0,x_2,0,x_4,0,\ldots)~:~x\in\ell_2\}.$$ Let $x=(x_1,0,x_3,0,\ldots) \in M$. Then $Tx=(x_3-x_1,x_1,x_5-x_3,\frac{x_3}{3},\ldots)$. Moreover, $T(M^\perp)=\{(x_2,x_2,x_4,x_4,\ldots)~:~ (0,x_2,0,x_4,\ldots)\in M^\perp \}$. Thus, $Tx+T(M^\perp)=\{(x_3+x_2-x_1,x_1+x_2,x_5+x_4-x_3,\frac{x_3}{3}+x_4,\ldots)~:~(0,x_2,0,x_4,\ldots)\in M^\perp \}$.

Let $z\in (T(x)+T(M^\perp))\cap M$. Then $z=Tx+Ty$ for some $y=(0,x_2,0,x_4,0,\ldots)\in M^\perp$. Therefore, $z=(x_3+x_2-x_1,x_1+x_2,x_5+x_4-x_3,\frac{x_3}{3}+x_4,\ldots)$. Since $z \in M$, we have $\frac{x_n}{n}+x_{n+1}=0$ for all odd integers $n$. Hence $x_{n+1}=-\frac{x_n}{n}$ for all odd integers $n$. So, $z=(x_3-2x_1,0,x_5-\frac{4x_3}{3},0,\ldots)$. Therefore $T(x)+T(M^\perp)$ is a singleton and $$T_{/M}(x_1,x_2,\ldots)=z=\Big(x_3-2x_1,0,x_5-\frac{4x_3}{3},0,\ldots\Big).$$
\end{example}

\noindent We proved that the `singleton condition' in the above theorem is a necessary and sufficient for the complementability of an operator $T=\begin{pmatrix}
A & B\\
C & D
\end{pmatrix}$ with $\mathcal R(D)$ closed.  Without the assumption of the closed range of the operator $D$, the `singleton condition' is only necessary but not sufficient, which is illustrated in the following example.

\begin{example}
	Let $T:\ell_2 \rightarrow \ell_2$ be defined by
	\[ T(x_1,x_2,\ldots) = \Big(x_1+\frac{x_2}{2}, x_1+\frac{x_2}{2},x_3+\frac{x_4}{4},x_3+\frac{x_4}{4},\ldots\Big). \]
Consider the subspaces	
	$$M=\{x=(x_1,0,x_3,0,\ldots)~:~x \in \ell_2\}$$ and 
	$$M^\perp =\{x=(0,x_2,0,x_4,0,\ldots)~:~x \in\ell_2 \}.$$ Therefore
		\begin{eqnarray*}
	 A(x_1,0,x_3,0, \ldots)&=&\Big(x_1,0,x_3,0, \ldots\Big);\\
	B(0,x_2,0,x_4,0, \ldots)&=&\Big(\frac{x_2}{2},0,\frac{x_4}{4},0, \ldots\Big);\\
	C(x_1,0,x_3,0, \ldots)&=&\Big(0,x_1,0,x_3,0, \ldots\Big);\\
	D(0,x_2,0,x_4,0, \ldots)&=&\Big(0,\frac{x_2}{2},0,\frac{x_4}{4},0, \ldots\Big).
		\end{eqnarray*}
	Clearly $\mathcal{R}(D)$ is not closed and $\mathcal{R}(C) \nsubseteq \mathcal{R}(D)$. So $T$ is not $M$-complementable.

	Let $x=(x_1,0,x_3,0,\ldots) \in M$. Then $Tx=(x_1,x_1,x_3,x_3,\ldots)$.  Thus $Tx+T(M^\perp)=\{(x_1+\frac{x_2}{2}, x_1+\frac{x_2}{2},x_3+\frac{x_4}{4},x_3+\frac{x_4}{4},\ldots)~:~(0,x_2,0,x_4,\ldots)\in M^\perp \}$.
	Suppose $z\in (T(x)+T(M^\perp))\cap M$. Then $z=Tx+Ty$ for some $y=(0,x_2,0,x_4,0,\ldots)\in M^\perp$. Therefore, $z=(x_1+\frac{x_2}{2}, x_1+\frac{x_2}{2},x_3+\frac{x_4}{4},x_3+\frac{x_4}{4},\ldots)$. Since $z \in M$, we have $\frac{x_{n+1}}{n+1}+x_{n}=0$ for all odd integers $n$. So, we get $z=0$. Therefore, $T(x)+T(M^\perp)$ is a singleton.

\end{example}
It is natural to ask, given $T\in \mathcal B(\mathcal H)$, whether there exist closed subspaces $M$ and $N$ such that $T$ is $(M,N)$ complementable. We affirmatively answer the question in the following results.
\begin{theorem}\label{excomp}
	Let  $T \in \mathcal{B}(\mathcal{H},\mathcal{K})$ and $N$ be a closed subspace of $\mathcal{K}$. Let $M=(T^{-1}(N\p))\p$. Then $T$ is $(M,N)$-complementable and $T_{/(M,N)}(x)=P_N Tx.$
\end{theorem}
\begin{proof}
	Let $T=\begin{pmatrix}
		A & B \\
		C & D
	\end{pmatrix}$. Given that $T(M\p)=N\p$, it follows that $\mathcal{R}(D)=\mathcal{R}(P_{N\p} T P_{M\p})=P_{N\p} T (M\p)=P_{N\p}(N\p)=N\p$, so $\mathcal{R}(D)$ is closed.
	For each $x\in M$, choose $y \in \{Tx + T(M\p) \} \cap N= \{Tx + N\p \} \cap N$. This implies that $y=Tx +z= P_N Tx + P_{N\p}Tx +z $,  for some  $z\in N\p$, so $y\in N$. As $z \in N\p$ and $y \in N$, we obtain $P_{N\p}Tx =-z$. Therefore, $y=P_N Tx + P_{N\p}Tx +z=P_N Tx$. Hence $\{Tx + T(M\p) \} \cap N =\{P_N Tx\}$, a singleton set. So by Theorem \ref{eqgm}, $T$ is $(M,N)$-complementable and  $T_{/(M,N)}(x)= P_N Tx$.
\end{proof}

\begin{corollary}
	Let  $T \in \mathcal{B}(\mathcal{H},\mathcal{K})$ and $M$ be a closed subspace of $\mathcal{H}$. Let $N=((T\s)^{-1}(M\p))\p$. Then $T$ is $(M,N)$-complementable.
\end{corollary}
\begin{proof}
	By Theorem \ref{excomp}, we have $T\s$ is $(N,M)$-complementable. Hence  $T$ is $(M,N)$-complementable by Theorem \ref{C}.
\end{proof}

\begin{proposition}
	Let  $T \in \mathcal{B}(\mathcal{H},\mathcal{K})$. Then there exist closed subspaces $M$ and $N$ such that $T$ is $(M,N)$-complementable.
\end{proposition}

\begin{proof}

	Let $x_1$ be a unit vector in $\mathcal{H}$ such that $T(x_1) \neq 0$. Consider an orthonormal basis $\{x_\lambda \}$ of $\mathcal{H}$ containing $x_1$, and let $\{y_\beta \}$ be an orthonormal basis of $\mathcal{K}$. Also, there exist non-zero coefficients $\beta _1$, $\beta _2$, $\ldots$, and vectors $y _1$, $y_2$, $\ldots$ in $\{y_\lambda \}$ such that $T(x_1)= \sum \limits _{m=1} ^{\infty} \beta_i y_i$. Define $M=(\text{span}\{x_1\})\p$ and $N=(\text{span}\{{y_1}\})\p$. The operator is expressed as $T=\begin{pmatrix}
		A & B\\
		C & D
	\end{pmatrix}$ with respect to the decomposition $M \oplus M\p$ and $N \oplus N\p$. Since $x_1 \in {M\p}$, we have $TP_{{M\p}}x_1=\sum \limits_{i=1}^{m}\beta_i y_i$. This implies $P_{{N\p}}TP_{{M\p}}x_1= \beta_1 y_1$, making $\mathcal{R}(P_{{N\p}}TP_{{M\p}})$ is an one-dimensional space. So, $\mathcal{R}(D)$ is a closed subspace.\\
	
	Let $x \in M$. Then $x=\sum \limits _{i=1}^{\infty}\langle x, x_i \rangle x_i$. Therefore, $Tx=\sum \limits _{i=1}^{\infty}\langle x, x_i \rangle Tx_i$. Since $Tx_i=\sum \limits _{i=1}^{\infty}\langle Tx_i, y_i \rangle y_i$, we have $$Tx= \sum \limits _{i=1}^{\infty}\big(\langle x, x_i \rangle \big(\sum \limits _{j=1}^{\infty}\langle Tx_i, y_j \rangle y_j \big) \big)=\sum \limits _{j=1}^{\infty}\bigg(\sum \limits _{i=1}^{\infty} \langle x, x_i \rangle \langle Tx_i, y_j \rangle \bigg) y_j.$$
	Let  $ w \in \{Tx + T(M\p) \} \cap N$. This implies that $w = Tx + Ty $ for some $y \in M\p$. Since $y \in M\p$, $y=\alpha _0 x_1$ for some $\alpha_0$. Therefore $$Ty=\alpha_0 T x_1=\alpha_0 \sum \limits _{k=1}^{\infty} \langle Tx_1, y_k \rangle y_k.$$ Consequently, $Tx+Ty= \sum \limits _{j=1}^{\infty}\bigg(\sum \limits _{i=1}^{\infty} \langle x, x_i \rangle \langle Tx_i, y_j \rangle \bigg) y_j + \alpha_0 \sum \limits _{k=1}^{\infty} \langle Tx_1, y_k \rangle y_k$. \\
	Since $Tx+Ty \in N$, the coefficient of $y_1=0$. Therefore $$\alpha_0 = -\frac{\sum \limits _{i=1}^{\infty} \langle x, x_i \rangle \langle Tx_i, y_1 \rangle}{\langle Tx_1, y_1 \rangle}.$$ This says that there is only one $y \in M\p$ such that $Tx+Ty \in N$. Therefore $\{Tx + T(M\p) \} \cap N =\{Tx+ T(\alpha_0 x_1)\}$, a singleton set.	Hence, by Theorem \ref{eqgm}, $T$ is $(M,N)$-complementable and \begin{align*}
		T_{/(M,N)}(x) &= Tx+T(\alpha_0 x_1)\\
		&=\sum \limits _{j=1}^{\infty}\bigg(\sum \limits _{i=1}^{\infty} \langle x, x_i \rangle \langle Tx_i, y_j \rangle \bigg) y_j + \alpha_0 Tx_1\\
		&=\sum \limits _{j=1}^{\infty}\bigg(\sum \limits _{i=1}^{\infty} \langle x, x_i \rangle \langle Tx_i, y_j \rangle \bigg) y_j - \Bigg(\frac{\sum \limits _{i=1}^{\infty} \langle x, x_i \rangle \langle Tx_i, y_1 \rangle}{\langle Tx_1, y_1 \rangle}\Bigg) Tx_1.
	\end{align*}
\noindent This completes the proof.
\end{proof}

		\begin{example}\label{ex6}
			Let $T:\ell_2 \to \ell_2$ be defined by  $$T(x_1,x_2,\ldots)=(x_1+x_2,x_1+x_2+x_3,x_2+x_3+x_4,x_3+x_4+x_5,\ldots).$$ 	and $M=N=(\{e_1\})\p$. With respect to the decompositions $M\oplus M^\perp$ and $N\oplus N^\perp$, we have the following operators :
		\begin{eqnarray*}
		A(0,x_2,x_3,\ldots)&=&(0,x_2+x_3,x_2+x_3+x_4,x_3+x_4+x_5,\ldots);\\  B(x_1,0,0,\ldots)&=&(0,x_1,0,0,\ldots);\\
		C(0,x_2,x_3,\ldots)&=&(x_2,0,0,\ldots);\\
		 D(x_1,0,0,\ldots)&=&(x_1,0,0,\ldots).
	\end{eqnarray*}
		\end{example}
		It is evident that $T$ is $(M,N)$-complementable, and all of the operators $A$, $B$, $C$, and $D$ are non-zero. This serves as a demonstration that the operators $A$, $B$, $C$, and $D$ in the previous proposition may not necessarily be zero.\\

		The current conceptualization and various characterizations of complementable operators are rooted in algebraic considerations, contingent upon the range of these operators. In this context, we introduce an analytical characterization for complementable operators. This characterization pivots on the properties of the image of the unit ball under these operators.
		
		 For a subspace $M$ of a Hilbert $\mathcal{H}$, we denote the set $\{x \in  M : \|x\| \leq 1 \}$ as $\mathcal{B}_M$.
		
		\begin{theorem}\label{char2}
			Let  $T \in \mathcal{B}(\mathcal{H},\mathcal{K})$. Then $T$ is $(M,N)$-complementable if and only if there exists $\lambda > 0$ such that  $C(\mathcal{B}_{M}) \subseteq \lambda D(\mathcal{B}_{ M\p})$ and $B\s(\mathcal{B}_{M}) \subseteq \lambda D\s(\mathcal{B}_{ M\p})$.
		\end{theorem}
		\begin{proof}
			Theorem \ref{dgls} says that  the inclusion $\mathcal{R}(C) \subseteq \mathcal{R}(D)$ gives the existence of an operator $Z \in \mathcal{B}(\mathcal{H})$ such that $C=DZ$. Suppose $C$ is non-zero. It follows that $Z$ is also non-zero, leading to the conclusion that $\|Z\| \neq 0$. Let $x \in \mathcal{B}_{\mathcal{H}}$. We then observe that $$Cx=DZx = \|Z\| \Big(D\Big(\frac{Zx}{\|Z\|}\Big)\Big).$$ Now we consider $y =\frac{Zx}{\|Z\|}$. So we have $\|y\| =\big \|\frac{Zx}{\|Z\|}\big \| \leq 1$, and $y \in \mathcal{R}(Z) \subseteq \mathcal{H}$. Therefore, we establish $C(\mathcal{B}_{\mathcal{H}}) \subseteq \|Z\| D(\mathcal{B}_{\mathcal{H}})$.
			
			Conversely, let us assume the existence of a positive constant $\lambda$ such that $C(\mathcal{B}_{\mathcal{H}}) \subseteq \lambda D(\mathcal{B}_{\mathcal{H}})$.
			Let $w \in \mathcal{R}(C)$, represented as $w=Cx$ for some $x \in \mathcal{H}$.  We observe that $C\Big(\frac{x}{\|x\|}\Big) \in C(\mathcal{B}_{\mathcal{H}}) \subseteq \lambda D(\mathcal{B}_{\mathcal{H}}) $. Accordingly, $ C\Big(\frac{x}{\|x\|}\Big) = \lambda D(y)$ for some $y \in \mathcal{B}_{\mathcal{H}}$. This implies $C{x} = \lambda \|x\| D(y)= D(\lambda \|x\| y)$. Hence we conclude $\mathcal{R}(C) \subseteq \mathcal{R}(D)$.
			
			Moreover,  $\| Z\| =\inf \{ \lambda >0 : C(\mathcal{B}_{\mathcal{H}}) \subseteq \lambda D(\mathcal{B}_{\mathcal{H}}) \}$. Indeed,  suppose that  $C(\mathcal{B}_{\mathcal{H}}) \subseteq \lambda D(\mathcal{B}_{\mathcal{H}})$ for some $\lambda >0$. For each $x \in \mathcal{B}_{\mathcal{H}}$, there exists $y\in \mathcal{B}_{\mathcal{H}} \cap \mathcal{N}(D)^\perp$ such that $Cx=\lambda Dy$ and also $Cx=DZx$ with $Zx \in \mathcal{N}(D)^\perp$. Therefore, we have $D(\lambda y-Zx)=0$, leading to $\lambda y-Zx \in \mathcal{N}(D) \cap \mathcal{N}(D)^\perp$. Therefore, $Zx=\lambda y$. As $\|y\| \leq 1$, we deduce $\|Zx\| = \lambda \|y\| \leq \lambda$ for all $x \in \mathcal{B}_{\mathcal{H}}$. Hence, we establish $\|Z\| \leq \lambda$.\\

			Expanding upon previous arguments, we establish the existence of positive constants $\lambda_1$ and $\lambda_2$ such that $C(\mathcal{B}_{M}) \subseteq \lambda_1 D(\mathcal{B}_{M\p})$ and $B\s(\mathcal{B}_{M}) \subseteq \lambda_2 D\s(\mathcal{B}_{M\p})$. Let $\lambda = \max \{\lambda_1, \lambda_2\} > 0$. Now, we have $\lambda_1 D(\mathcal{B}_{M\p}) \subseteq \lambda D(\mathcal{B}_{M\p})$ and $\lambda_2 D\s(\mathcal{B}_{M\p}) \subseteq \lambda D\s(\mathcal{B}_{M\p})$. Hence we have $C(\mathcal{B}_{M}) \subseteq \lambda D(\mathcal{B}_{M\p})$ and $B\s(\mathcal{B}_{M}) \subseteq \lambda D\s(\mathcal{B}_{M\p})$.
		\end{proof}

		 	Let $T\in \mathcal{B}(\mathcal{H},\mathcal{K})$ be an operator. The reduced minimum modulus of $T$ is defined as $$\gamma (T)=\inf \{\|Tx\| : x\in \mathcal{H} \cap \mathcal{N}(T)\p \}.$$ 
Moreover, $\mathcal{R}(T) $ is closed if and only if $\gamma (T)>0$ \cite{minmod}. We now characterize algebraically $(M,N)$-complementable  operator 
$T=\begin{pmatrix}
A & B\\
C & D
\end{pmatrix}$ with respect to the reduced minimum modulus of the operator $D$.		    
		\begin{theorem}\label{thm34}
			Let  $T \in \mathcal{B}(\mathcal{H},\mathcal{K})$ with $\mathcal{R}(D)$ closed. Then $T$ is $(M,N)$-complementable if and only if $C(\mathcal{B}_{M}) \subseteq \frac{\|C\|}{\gamma (D)}D(\mathcal{B}_{ M\p})$ and $B\s(\mathcal{B}_{M}) \subseteq \frac{\|B\|}{\gamma (D\s)}D\s(\mathcal{B}_{M\p})$.
		\end{theorem}

		\begin{proof}
			Let $T$ be $(M,N)$-complementable and $x \in
			\mathcal{B}_\mathcal H$.  Since			 $\mathcal{R}(C) \subseteq \mathcal{R}(D)$, there exists $y \in (\mathcal{N}(D)\p) \cap M\p$ such that $Cx = Dy$. As $\mathcal{R}(D)$ is closed, we have $$\|y\| \gamma (D) \leq \|y\| \|D \big(\frac{y}{\|y\|}\big)\|=\|Dy\| =\|Cx\| \leq \|C\|.$$ Therefore the norm of $y$ is bounded by $\frac{\|C\|}{\gamma(D)}$. Therefore, $y$ can be expressed as $y = k y_0$, where $0 \leq k \leq \frac{\|C\|}{\gamma(D)}$, and $y_0 = \frac{y}{\|y\|}$. Thus, $y = \frac{\|C\|}{\gamma(D)}\left(k \frac{\gamma(D)}{\|C\|}y_0\right)$. Since $k \frac{\gamma(D)}{\|C\|} \leq 1$, it follows that $y = \frac{\|C\|}{\gamma(D)}z$ for some $\|z\| \leq 1$. Hence, $C(\mathcal{B}_{M})$ is contained in $\frac{\|C\|}{\gamma(D)}D(\mathcal{B}_{M^{\perp}})$. 		In a similar way, the same argument can be made for $B^*$ and $D^*$, leading to the inclusion $$B^*(\mathcal{B}_{M}) \subseteq \frac{\|B^*\|}{\gamma(D^*)}D^*(\mathcal{B}_{ M^{\perp}}).$$
			
			Conversely, if $C(\mathcal{B}_{M}) \subseteq \frac{\|C\|}{\gamma(D)}D(\mathcal{B}_{M^{\perp}})$ and $B^*(\mathcal{B}_{M}) \subseteq \frac{\|B^*\|}{\gamma(D^*)}D^*(\mathcal{B}_{ M^{\perp}})$, then by Theorem \ref{char2}, $T$ is $(M,N)$-complementable.
		\end{proof}

		\section{\textbf{Schur Complements of  EP and  Hypo-EP Operators}}
\noindent		Motivated by the Schur determinant formula
		$$det(T)=det(A) \ det (T_{/A})$$ for matrices, where $T_{/A}$ is given in (\ref{equation-1}), we analyze the  closedness of the range of  $(M, N)$-complementable operator.
		\begin{theorem}\label{t312}
			Let  $T \in \mathcal{B}(\mathcal{H},\mathcal{K})$ be $(M,N)$-complementable. Then $\mathcal{R}(T)$ is closed if and only if $\mathcal{R}(T_{/(M,N)})$ and $\mathcal{R}(D)$ are closed.
		\end{theorem}
		\begin{proof}
			Let $Z$ and $Y$ be such that $C=DZ$ and $B=YD$. This implies the following matrix equality:
			\begin{equation*}
				\begin{pmatrix}
					A & B\\
					C & D
				\end{pmatrix} =
				\begin{pmatrix}
					I & Y\\
					0 & I
				\end{pmatrix}
				\begin{pmatrix}
					T_{/(M,N)} & 0\\
					0 & D
				\end{pmatrix}
				\begin{pmatrix}
					I & 0\\
					Z & I
				\end{pmatrix}
			\end{equation*}
			where $\begin{pmatrix}
				I & Y\\
				0 & I
			\end{pmatrix}$ and $\begin{pmatrix}
				I & 0\\
				Z & I
			\end{pmatrix}$ are invertible.
			
			\noindent Assuming $\mathcal{R}(T)$ is closed, we can express the above equality as:
			\begin{equation*}
				\begin{split}
					\begin{pmatrix}
						T_{/(M,N)} & 0\\
						0 & D
					\end{pmatrix} &=
					\begin{pmatrix}
						I & -Y\\
						0 & I
					\end{pmatrix}
					\begin{pmatrix}
						A & B\\
						C & D
					\end{pmatrix}
					\begin{pmatrix}
						I & 0\\
						-Z & I
					\end{pmatrix}.
				\end{split}
			\end{equation*}
			Let $P=\begin{pmatrix}
				I & Y\\
				0 & I
			\end{pmatrix}$, $Q=\begin{pmatrix}
				I & 0\\
				Z & I
			\end{pmatrix}$ and $U=\begin{pmatrix}
				T_{/(M,N)} & 0\\
				0 & D
			\end{pmatrix}$.
			
			\noindent Now, analyzing the range of $U$, we have
			\begin{equation*}
				\begin{split}
					\mathcal{R}(U) &= U(\mathcal{H}) \\
					&= (P^{-1}TQ^{-1}) (\mathcal{H}) \\
					&= P^{-1}T(Q^{-1} (\mathcal{H})) \\
					&= (P^{-1}T) (\mathcal{H}) \\
					&= P^{-1} (\mathcal R(T)).
				\end{split}
			\end{equation*}
			Since $\mathcal{R}(T)$ is closed and  $P$ is bounded,  $\mathcal{R}(U)$ is closed. However, $\mathcal{R}(U)=\mathcal{R}(T_{/(M,N)}) \oplus \mathcal{R}(D)$. Therefore, both $\mathcal{R}(T_{/(M,N)})$ and $\mathcal{R}(D)$ are closed.
			
			Conversely, assuming $\mathcal{R}(T_{/(M,N)})$ and $\mathcal{R}(D)$ are closed, we observe that $\mathcal{R}(T)= PUQ(\mathcal{H})=PU(\mathcal{H})$. Now $U(\mathcal{H})= \mathcal{R}(T_{/(M,N)}) \oplus \mathcal{R}(D)$. Given the closure of $\mathcal{R}(T_{/(M,N)})$ and $ \mathcal{R}(D)$, $U(\mathcal{H})= \mathcal{R}(T_{/(M,N)}) \oplus \mathcal{R}(D)$ is closed. Because $P^{-1}$ is bounded, $P(U(\mathcal{H}))$ is closed. Therefore, $\mathcal{R}(T)$ is closed.
		\end{proof}
		
The following example provides a $(M,N)$-complementable operator whose range is not closed.  \begin{example}\label{ex1}
			Let $T:\ell_2 \to \ell_2$ be defined by $$T(x_1,x_2,\ldots)= \Big(x_1,x_2,\frac{x_3}{3},x_4,\frac{x_5}{5},x_6,\ldots\Big).$$ Consider the subspaces $$M=\{(x_1,x_1,x_3,x_3,\ldots)\} \text { and } M\p =\{(x_2,-x_2,x_4,-x_4,\ldots)\}.$$ We have $T=\begin{pmatrix}
			A & B\\
			C & D
			\end{pmatrix}$
			where  $A$, $B$, $C$ and $D$ are given by
			\begin{eqnarray*}
				A(x_1,x_1,x_3,x_3\ldots)&=&\frac{1}{2}\Big(2x_1,2x_1,\frac{4}{3}x_3,\frac{4}{3}x_3,\ldots\Big); \\ B(x_1,-x_1,x_3,-x_3\ldots)&=&\frac{1}{2}\Big(0,0,-\frac{2}{3}x_3,-\frac{2}{3}x_3,-\frac{4}{5}x_5,-\frac{4}{5}x_5,\ldots\Big);\\ C(x_1,x_1,x_3,x_3\ldots)&=&\frac{1}{2}\Big(0,0,-\frac{2}{3}x_3,\frac{2}{3}x_3,-\frac{4}{5}x_5,\frac{4}{5}x_5,\ldots\Big); \\ D(x_1,-x_1,x_3,-x_3\ldots)&=&\frac{1}{2}\Big(2x_1,-2x_1,\frac{4}{3}x_3,-\frac{4}{3}x_3,\ldots\Big).
			\end{eqnarray*}
		\end{example}
		Clearly, $\mathcal{R}(T)$ is not closed and $\mathcal{R}(D)$ is closed. Now by definition we have $$T_{/M}(x_1,x_2,x_3,x_4,\ldots)=\frac{1}{2}\Big(0,0,x_3,x_3,\frac{4}{6}x_5,\frac{4}{6}x_5,\frac{4}{8}x_7,\frac{4}{8}x_7,\ldots \Big).$$ So $\mathcal{R}(T_{/M})$ is not closed.

  It is established in  \cite{Meenakshi} that under certain conditions a Schur complement of an $EP$-matrix is $EP$ again.  Motivated by this work, we now discuss Schur complements of $EP$ and hypo-$EP$ operators. 
  An operator			$T\in \mathcal{B({H})}$ is said to be $EP$ (hypo-$EP$) if $\mathcal{R}(T)$ is closed and  $\mathcal{R}(T)=\mathcal{R}(T\s)$ ($\mathcal{R}(T)\subseteq \mathcal{R}(T\s)$) \cite{Campbell,Itoh}.

		\begin{theorem}\label{ep}
			Let $T \in \mathcal{B({H})}$ be a $(M,N)$-complementable operator. 
			If $T$ is hypo-$EP$, then $T_{/(M,N)}$ is hypo-$EP$. \\
			Moreover, if $\mathcal{R}(C\s) \subseteq \mathcal{R}((T_{/(M,N)})\s)$, then the following are equivalent:
			\begin{enumerate}
				\item $T$ is hypo-$EP$;
				\item $T_{/(M,N)}$ and $D$ are hypo-$EP$ and $\mathcal{R}(B) \subseteq \mathcal{R}(T_{/(M,N)})$;
				\item $\begin{pmatrix}
					T_{/(M,N)} & 0\\
					C & D
				\end{pmatrix} $ and $\begin{pmatrix}
					T_{/(M,N)} & B\\
					0 & D
				\end{pmatrix} $ are hypo-$EP$.
			\end{enumerate}
		\end{theorem}
		\begin{proof}
			Suppose that $T$ is hypo-$EP$.  Theorem \ref{C} establishes that $\mathcal{R}({T_{/(M,N)}})=\mathcal{R}(T) \cap N$ and $\mathcal{R}((T\s_{/(M,N)}))=\mathcal{R}(({T_{/(M,N)}})\s)=\mathcal{R}(T\s) \cap N$. The hypo-$EP$ property of $T$ implies $\mathcal{R}({T_{/(M,N)}})\subseteq \mathcal{R}(({T_{/(M,N)}})\s)$. Furthermore, given the closed nature of both $\mathcal{R}(T)$ and $N$, it follows that $\mathcal{R}({T_{/(M,N)}})$ is closed. Therefore, ${T_{/(M,N)}}$ is hypo-$EP$.\\
			
			$(1) \implies (2)$: Already we have proved that ${T_{/(M,N)}}$ is hypo-$EP$, if $T$ is hypo-$EP$.  Considering the matrix operators $P=\begin{pmatrix}
				I & BD\+ \\
				0 & I
			\end{pmatrix}$, $Q=\begin{pmatrix}
				I & 0\\
				C{T_{/(M,N)}}\+ & I
			\end{pmatrix}$, and $L=\begin{pmatrix}
				{T_{/(M,N)}} & 0\\
				0 & D
			\end{pmatrix}$, we can express  $P QL$ as
			
			\begin{align*}
				P QL&=\begin{pmatrix}
					I+BD\+C({T_{/(M,N)}})\+ & BD\+\\
					C({T_{/(M,N)}})\+ & I
				\end{pmatrix}\begin{pmatrix}
					{T_{/(M,N)}} & 0\\
					0 & D
				\end{pmatrix}\\
				&=\begin{pmatrix}
					{T_{/(M,N)}}+BD\+C({T_{/(M,N)}})\+{T_{/(M,N)}} & BD\+D\\
					C{T_{/(M,N)}}({T_{/(M,N)}})\+ & D
				\end{pmatrix}\\
				&=\begin{pmatrix}
					A & B\\
					C & D
				\end{pmatrix}.
			\end{align*}
			
			This implies that $\mathcal{N}(L)=\mathcal{N}(T)$. Taking the orthogonal complement on both sides yields $\mathcal{R}(T\s)=\mathcal{R}(L\s)$. Since $T$ is hypo-$EP$, we have $\mathcal{R}(T)\subseteq \mathcal{R}(L\s)$. Additionally, $\mathcal{R}(L)$ is closed. Thus, $$T\s=T\s L\+L=\begin{pmatrix}$
				$A\s ({T_{/(M,N)}})\+ {T_{/(M,N)}} & C\s D\+ D\\
				B\s ({T_{/(M,N)}})\+ {T_{/(M,N)}} & D\s D\+ D
			\end{pmatrix}.$$
						This leads to the relationships $B\s {=} B\s ({T_{/(M,N)}})\+ {T_{/(M,N)}}$ and $D\s {=} D\s D\+ D$. As a result, $ \mathcal{R}(B) \subseteq \mathcal{R}({T_{/(M,N)}})$ and $\mathcal{R}(D) \subseteq \mathcal{R}(D\s)$.
			
			$(2) \implies (1)$: Given that $\mathcal{R}({T_{/(M,N)}})$ and $\mathcal{R}(D)$ are closed, as established by Theorem \ref{t312}, it follows that $\mathcal{R}(T)$ is also closed. By applying Theorem \ref{t215}, we express $T\+$ as
						\[ T\+=\begin{pmatrix}
				({T_{/(M,N)}})\+ & -({T_{/(M,N)}})\+ BD\+\\
				-D\+C({T_{/(M,N)}})\+ & D\+ +D\+ C({T_{/(M,N)}})\+BD\+
			\end{pmatrix} .\]
			
			Utilizing the conditions $\mathcal{R}(C) \subseteq \mathcal{R}(D)$, $\mathcal{R}(B\s) \subseteq \mathcal{R}(D\s)$, $\mathcal{R}(C\s) \subseteq \mathcal{R}({T_{/(M,N)}}\s)$ and $\mathcal{R}(B) \subseteq \mathcal{R}({T_{/(M,N)}})$, we deduce 
			$$TT\+=\begin{pmatrix}
				{T_{/(M,N)}}({T_{/(M,N)}})\+ & 0\\
				0 & DD\+
			\end{pmatrix} \ \text{ and } \ 
			T\+T=\begin{pmatrix}
				({T_{/(M,N)}})\+ {T_{/(M,N)}} & 0\\
				0 & D\+D
			\end{pmatrix}.$$ Given that ${T_{/(M,N)}}$ and $D$ are hypo-$EP$, it follows that $\mathcal{R}(TT\+) \subseteq \mathcal{R}(T\+ T)$. Accordingly, $ \mathcal{R}(T) \subseteq \mathcal{R}(T\s)$.  Thus $ T$ is hypo-$EP$.
			
			$(2) \iff (3)$: Given that $T$ is $(M,N)$-complementable. The first implication indicates that $\begin{pmatrix} {T_{/(M,N)}} & 0\\ C & D \end{pmatrix}$ and $\begin{pmatrix} {T_{/(M,N)}} & B\\ 0 & D \end{pmatrix}$ are also $(M,N)$-complementable. Both $\begin{pmatrix} {T_{/(M,N)}} & 0\\ C & D \end{pmatrix}$ and $\begin{pmatrix} {T_{/(M,N)}} & B\\ 0 & D \end{pmatrix}$ satisfy all the conditions established for $T$ in the hypothesis. Additionally, $({T_{/(M,N)}})_{/(M,N)}=T_{/(M,N)}$. Thus, by applying 	$(1) \iff (2)$ for operators $\begin{pmatrix} T_{/(M,N)} & 0\\ C & D \end{pmatrix}$ and $\begin{pmatrix} {T_{/(M,N)}} & B\\ 0 & D \end{pmatrix}$, we obtain that $\begin{pmatrix} {T_{/(M,N)}} & 0\\ C & D \end{pmatrix}$ and $\begin{pmatrix} {T_{/(M,N)}} & B\\ 0 & D \end{pmatrix}$ are hypo-$EP$.

		\end{proof}
		
		\begin{corollary}
			Let $T \in \mathcal{B({H})}$ be a $(M,N)$-complementable operator. 
			If $T$ is $EP$, then $T_{/(M,N)}$ is $EP$.\\
			Moreover, if $\mathcal{R}(C\s) \subseteq \mathcal{R}((T_{/(M,N)})\s)$ then the following are equivalent:
			\begin{enumerate}
				\item $T$ is $EP$ ;
				\item $T_{/(M,N)}$ and $D$ are $EP$ and $\mathcal{R}(B)  \subseteq \mathcal{R}(T_{/(M,N)})$ ;
				\item $\begin{pmatrix}
					T_{/(M,N)} & 0\\
					C & D
				\end{pmatrix} $ and $\begin{pmatrix}
					T_{/(M,N)} & B\\
					0 & D
				\end{pmatrix} $ are $EP$.
			\end{enumerate}
		\end{corollary}
		\begin{proof}
			By Theorem \ref{ep}, it is easy to see that $T_{/(M,N)}$ is $EP$. \\
			$(1) {\iff} (2)$: It is easy to see that $T$ is EP if and only if both $T$ and $T^*$ are hypo-EP.  Hence by Theorem \ref{ep}, $T$ is EP if and only if ${T_{/(M,N)}}$ and  $D$ are EP, and hence $\mathcal{R}(B) \subseteq \mathcal{R}({T_{/(M,N)}})$.
			
			$(2) {\iff} (3)$: Both $\begin{pmatrix}
				T_{/(M,N)} & 0\\
				C & D
			\end{pmatrix}$ and $\begin{pmatrix}
				T_{/(M,N)} & B\\
				0 & D
			\end{pmatrix}$ are $(M, N)$-complementable, and $(T_{/(M,N)})_{/(M,N)}=T_{/(M,N)}$. 
			
			Applying $(1) {\iff} (2)$ for operators $\begin{pmatrix}
				T_{/(M,N)} & 0\\
				C & D
			\end{pmatrix}$ and $\begin{pmatrix}
				T_{/(M,N)} & B\\
				0 & D
			\end{pmatrix}$, we conclude that  $\begin{pmatrix}
				T_{/(M,N)} & 0\\
				C & D
			\end{pmatrix}$ and $\begin{pmatrix}
				T_{/(M,N)} & B\\
				0 & D
			\end{pmatrix}$ are $EP$.

		\end{proof}

		Now we illustrate Theorem \ref{ep} with the following example.  
		
		 \begin{example}
			Let $T:\ell_2 \oplus \ell_2 \to \ell_2 \oplus \ell_2$ be defined by $$ T((x_1, x_2,\ldots) \oplus (y_1,y_2,\ldots ))=\Big( \Big(x_1, x_1+\frac{x_2}{2}+y_1,x_2+\frac{2x_3}{3}+y_2,\ldots\Big) \oplus \Big(0,x_1+y_1,x_2+y_2,\ldots \Big) \Big).$$
			 Consider the subspaces $$M=\ell_2 \oplus 0 \text { and } M\p =0 \oplus \ell_2.$$ We have $T=\begin{pmatrix}
				A & B\\
				C & D
			\end{pmatrix}$
			where  $A$, $B$, $C$ and $D$ are given by
			\begin{eqnarray*}
				A((x_1,x_2,\ldots)\oplus 0)&=&\bigg(\bigg(x_1,x_1+\frac{x_2}{2},x_2+\frac{2x_3}{3},\ldots\bigg) \oplus 0 \bigg);\\
			 B(0 \oplus (y_1,y_2,\ldots))&=&\big( (0,y_1,y_2,\ldots) \oplus 0 \big);\\
			  C((x_1,x_2,\ldots)\oplus 0)&=& \big(0 \oplus (0,x_1,x_2,\ldots) \big);  \\
			   D(0 \oplus (y_1,y_2,\ldots))&=&\big(0 \oplus (0,y_1,y_2,\ldots)\big).
			\end{eqnarray*}
		\end{example}

		Clearly, $T$ is $M$-complementable and $${T_{/M}}\big((x_1,x_2,x_3,\ldots) \oplus 0 \big)=\bigg(\bigg(x_1,\frac{x_2}{2},\frac{2x_3}{3}\ldots \bigg) \oplus 0 \bigg).$$ Also, ${T_{/M}}\s={T_{/M}}$, $\mathcal{R}(D)$ is closed and $\mathcal{R}(C\s) \subseteq \mathcal{R}({T_{/M}}\s)$. Now let $x \oplus y \in \mathcal{N}(T)$, then $x=y=0$. Then $0=\mathcal{N}(T)\subseteq \mathcal{N}(T\s)$.  Now ${T_{/M}}:\ell_2 \to \ell_2$ is invertible, so $\mathcal{R}({T_{/M}})$ is closed. Also, $0=\mathcal{N}({T_{/M}})\subseteq \mathcal{N}({T_{/M}}\s)$. Also by Theorem \ref{t312}, $\mathcal{R}(T)$, $\mathcal{R}(T_1)$, and $\mathcal{R}(T_2)$ are closed, where $T_1=\begin{pmatrix}
			{T_{/M}}&0\\ C&D
		\end{pmatrix}$ and  $T_2=\begin{pmatrix}
			{T_{/M}}&B\\ 0 &D
		\end{pmatrix}.$ 
		So we have $T$, ${T_{/M}}$, $D$, $T_1$ and $T_2$ are hypo-$EP$ and $\mathcal{R}(B)\subseteq \mathcal{R}({T_{/M}})$.

		The next example shows that the condition $\mathcal{R}(C\s)\nsubseteq \mathcal{R}((T_{/M})\s)$ is necessary for the implication 	$(1)$ if and only if $(2)$ in Theorem \ref{ep}.
		
		 \begin{example}\label{ex3}
			Let $T:\ell_2 \to \ell_2$ be defined by $$ T(x_1,x_2,\ldots)=(x_1+x_2+x_3+x_4,0,x_2+x_3+x_4,x_2+x_3+x_4,x_5,x_6,\ldots).$$
			 Consider the subspaces $$M=\{(x_1,x_2,0,0,\ldots)\} \text { and } M\p =\{(0,0,x_3,x_4,\ldots)\}.$$ We have $T=\begin{pmatrix}
				A & B\\
				C & D
			\end{pmatrix}$
			where  $A$, $B$, $C$ and $D$ are given by
			\begin{eqnarray*}
				A(x_1,x_2,0,0,\ldots)&=&(x_1+x_2,0,0,0 \ldots);\\
			B(0,0,x_3,x_4,\ldots)&=&(x_3+x_4,0,0,0 \ldots);\\
				C(x_1,x_2,0,0,\ldots)&=&(0,0,x_2,x_2,0,0,\ldots);\\  
				D(0,0,x_3,x_4,\ldots)&=&(0,0,x_3+x_4,x_3+x_4,x_5,x_6\ldots).\\
			\end{eqnarray*}
		\end{example}
		Clearly, $T$ is $M$-complementable and $$T_{/M}(x_1,x_2,0,0,\ldots)=(x_1,0,0,\ldots).$$ But $\mathcal{R}(C\s)\nsubseteq \mathcal{R}((T_{/M})\s)$. Now, $D$ and $T_{/M}$ are hypo-$EP$ ($EP$), but $T$ is not hypo-$EP$ ($EP$). This shows that the range space conditions are necessary for Theorem \ref{ep}.

	 	The next example shows that the condition $\mathcal{R}(C\s)\nsubseteq \mathcal{R}((T_{/M})\s)$ is necessary for the implication 	$(2)$ if and only if $(3)$ in Theorem \ref{ep}.
	 	
	 \begin{example}\label{ex4}
	 	Let $T:\ell_2 \to \ell_2$ be defined by $$ T(x_1,x_2,\ldots)=(x_1+x_2+x_3+x_4,x_1+x_3+x_4,x_1+x_2+x_3+x_4,x_1+x_2+x_3+x_4,x_5,x_6,\ldots).$$ with respect to the subspaces $$M=\{(x_1,x_2,0,0,\ldots)\} \text { and } M\p =\{(0,0,x_3,x_4,\ldots)\}.$$ We have $T=\begin{pmatrix}
	 		A & B\\
	 		C & D
	 	\end{pmatrix}$
	 	where  $A$, $B$, $C$ and $D$ are given by

	 	\begin{eqnarray*}
	 		A(x_1,x_2,0,0,\ldots)&=&(x_1+x_2,x_1,0,0,\ldots); \\ B(0,0,x_3,x_4,\ldots)&=&(x_3+x_4,x_3+x_4,0,0,\ldots);\\ C(x_1,x_2,0,0,\ldots)&=&(0,0,x_1+x_2,x_1+x_2,0,0,\ldots);\\ D(0,0,x_3,x_4,\ldots)&=&(0,0,x_4+x_4,x_3+x_4,x_5,x_6\ldots).\\
	 	\end{eqnarray*}
	 \end{example}
	 Clearly, $T$ is $M$-complementable and $$T_{/M}(x_1,x_2,0,0,\ldots)=(0,-x_2,0,\ldots).$$ But $\mathcal{R}(B)\nsubseteq \mathcal{R}(T_{/M})$ and $\mathcal{R}(C\s)\nsubseteq \mathcal{R}((T_{/M})\s)$. Clearly, $D$, $T_{/M}$ and  $T$ are hypo-$EP$ ($EP$), but $\begin{pmatrix}
	 	T_{/M}&0\\ C&D
	 \end{pmatrix}$ and $\begin{pmatrix}
	 	T_{/M}&B\\ 0&D
	 \end{pmatrix}$ are not hypo-$EP$ ($EP$).

	\begin{center}
		{\bf Acknowledgments }
	\end{center}
	The first author thanks the National Institute of Technology Karnataka (NITK), Surathkal for the financial support. The present work of the second author was partially supported by Science and Engineering Research Board (SERB), Department of Science and Technology, Government of India (Reference Number: MTR/2023/000471) under the scheme ``Mathematical Research Impact Centric Support (MATRICS)".
\end{sloppypar}


\begin{thebibliography}{99}
		\bibitem{Anderson}
		W. N. Anderson, Jr.,
		\newblock Shorted operators.
		\newblock {\em SIAM J. Appl. Math.}, 20:520--525, 1971.
		
		
		\bibitem{Ando}
		T. Ando.,
		\newblock Generalized {S}chur complements
		\newblock  {\em Linear Algebra Appl.},  27:173--186, 1979.
		
		
		
		\bibitem{Antezana}
		J. Antezana, G. Corach and D. Stojanoff.,
		\newblock Bilateral shorted operators and parallel sums.
		\newblock {\em Linear Algebra Appl.}, 414:570--588, 2006.
		
		
		
		\bibitem{Arias}
		M. L. Arias, G. Corach and A. Maestripieri.,
		\newblock On complementable operators in the sense of {T}. {A}ndo.
		\newblock {\em Linear Algebra Appl.}, 594:158--176, 2020.
	
		\bibitem{Albert1}
A. Arthur.,
\newblock Conditions for positive and nonnegative definiteness in terms
of pseudoinverses.
\newblock {\em SIAM J. Appl. Math.}, 17:434--440,
1969.


\bibitem{Butler}
C. A. Butler and T.D. Morley.,
\newblock Six generalized {S}chur complements
\newblock {\em Linear Algebra Appl.}, 106:259--269, 1988.


\bibitem{Campbell}
S. L. Campbell and C. D. Meyer.,
\newblock  {$EP$} operators and generalized inverses. {\em Canad. Math. Bull.}, 18(3):327--333, 1975.




\bibitem{Carlson1}
D. Carlson.,
\newblock What are {S}chur complements, anyway?.
\newblock {\em Linear Algebra Appl.}, 74:257--275, 1986.


\bibitem{Corach2}
G. Corach, A. Maestripieri and D. Stojanoff.,
\newblock {Generalized {S}chur complements and oblique projections.}
\newblock {\em Linear Algebra Appl.},  341, 259--272, 2002.

\bibitem{Corach1}
G. Corach, A. Maestripieri and D. Stojanoff.,
\newblock {Oblique projections and {S}chur complements}.
\newblock {\em Acta Sci. Math. (Szeged)},  67, 337--356, 2001.



\bibitem{Cottle}
R. W. Cottle.,
\newblock {On manifestations of the {S}chur complement.}
\newblock {\em Rend. Sem. Mat. Fis. Milano },  45, 31--40, 1975.



\bibitem{Gantmacher}
F. R. Gantmacher.,	
\newblock{\em The theory of matrices. {V}ol. 1}.
\newblock {AMS Chelsea Publishing, Providence, RI, 1998}.

\bibitem{Giribet}
J. I. Giribet, A. Maestripieri and F. M.
Per\'{\i}a.,
\newblock {Shorting selfadjoint operators in {H}ilbert spaces.}
\newblock {\em Linear Algebra Appl.},  428, 1899--1911, 2008.





\bibitem{MR3314341}	
N. C. Gonz\'{a}lez, M. F. Mart\'{\i}nez-Serrano and J. Robles.,
\newblock Expressions for the {M}oore-{P}enrose inverse of block
matrices involving the {S}chur complement.
\newblock {\em Linear Algebra Appl.}, 471:353--368,  2015.

\bibitem{Emilie}
E. V. Haynsworth.,
\newblock Determination of the inertia of a partitioned {H}ermitian
matrix
\newblock {\em Linear Algebra Appl.}, 1:73--81, 1968.

\bibitem{minmod}		
I. S. Hwang and W. Y. Lee.,
\newblock {The reduced minimum modulus of operators.}
\newblock {\em J. Math. Anal. Appl.},  267, 679--694, 2002.	






\bibitem{Itoh}
M. Itoh.,
\newblock  On some {EP} operators.
{\em Canad. Math. Bull.}, 16(1):49--56, 2005.

\bibitem{Krein}
M. Krein.,
\newblock The theory of self-adjoint extensions of semi-bounded
{H}ermitian transformations and its applications {I}.
\newblock {\em Rec. Math. [Mat. Sbornik] N.S.}, 20(62):431--495, 1947.




\bibitem{Mitra}	
S. K. Mitra and  M. L. Puri.,
\newblock	Shorted matrices---an extended concept and some applications.
\newblock {\em Linear Algebra Appl.}, 42:57--79,  1982.	




\bibitem{Meenakshi}
Ar. Meenakshi.,
\newblock On {S}chur complements in an {EP} matrix.
\newblock {\em Period. Math. Hungar.}, 16:193--200, 1985.


\bibitem{Ouellette}
D. V. Ouellette.,
\newblock {Schur complements and statistics.}
\newblock {\em  Linear Algebra Appl.},36, 187-295,1981.	



\bibitem{Pekarev}
E. L. Pekarev., 
\newblock {Shorts of operators and some extremal problems.}
\newblock {\em Acta Sci. Math. (Szeged) },  56:1-2,  147-163, 1992.


\bibitem{Ptak}
V. Ptak.,
\newblock {Extremal operators and oblique projections.}
\newblock {\em Casopis pro pestovani Matematiky },  110,  343--350, 1985.




\bibitem{Schur}
J. Schur.,
\newblock \"{U}ber {P}otenzreihen, die im {I}nnern des {E}inheitskreises
beschr\"{a}nkt sind.
\newblock {\em J. Reine Angew. Math.}, 147:205--232, 1917.	


\bibitem{sylvester}		
J. J. Sylvester.,
\newblock {On the relation between the minor determinants of linearly equivalent quadratic functions.}
\newblock {\em The London, Edinburgh, and Dublin Philosophical Magazine and Journal of Science},  1(4), 295–305, 1851.	




\bibitem{Zhang}		
F. Zhang.,
\newblock {The Schur Complement and Its Applications.}
\newblock {\em Springer New York, NY}, 1571-5698.















\end{thebibliography}
\end{document}